\documentclass{conm-p-l}
\usepackage{mathrsfs}
\usepackage{amsmath, amsthm, amssymb, enumerate, mathtools}

\usepackage{color}

\numberwithin{equation}{section}

\title[Con\-for\-mal parity]{Con\-for\-mal parity}

\author[A.\,Derdzinski, S.\,Kim, J.H.\,Park]{Andrzej Derdzinski$^{1}$,
Sinhwi Kim$^{2}$, and JeongHyeong Park$^{2}$}
\address{$^{1}$\nnh\ Department\nnh\ of\nnh\ Mathematics,\nnh\ The\nnh\
Ohio\nnh\ State\nnh\ University,\nnh\ Columbus,\nnh\ OH\nnh\ 43210,\nnh\ USA}
\email{andrzej@math.ohio-state.edu}

\address{$^{2}$ Department of Mathematics, Sungkyunkwan University, Suwon, 16419, Korea}
\email{kimsinhwi@skku.edu, parkj@skku.edu}


\def\bbC{{\mathchoice {\setbox0=\hbox{$\displaystyle\rm C$}\hbox{\hbox
to0pt{\kern0.4\wd0\vrule height0.9\ht0\hss}\box0}}
{\setbox0=\hbox{$\textstyle\rm C$}\hbox{\hbox
to0pt{\kern0.4\wd0\vrule height0.9\ht0\hss}\box0}}
{\setbox0=\hbox{$\scriptstyle\rm C$}\hbox{\hbox
to0pt{\kern0.4\wd0\vrule height0.9\ht0\hss}\box0}}
{\setbox0=\hbox{$\scriptscriptstyle\rm C$}\hbox{\hbox
to0pt{\kern0.4\wd0\vrule height0.9\ht0\hss}\box0}}}}

\def\bbR{\mathrm{I\!R}}

\def\rtr{\bbR\nh^3}

\def\txm{{T\hskip-2.9pt_x\w\hn M}}

\newcommand{\sca}{\mathrm{s}}
\newcommand{\w}{^{\phantom i}}
\newcommand{\nnh}{\hskip-1pt}

\def\ya{a}
\def\yb{b}
\def\yc{c}

\def\x{x}
\def\y{y}
\def\z{z}
\def\t{t}

\def\hj{H}

\def\n{\nabla\hskip-2.4pt}

\def\ly{\lambda}
\def\sy{\sigma}
\def\my{\mu}

\def\ky{\kappa}

\def\hi{\varphi}

\def\bde{\overline{\nnh\Delta\nnh}}

\usepackage{color}

\newcommand{\ric}{\mathrm{r}}

\def\cwedge{\bigcirc\kern-1.07em\wedge\ }
\newcommand{\hyp}{\hskip.5pt\vbox
{\hbox{\vrule width2.5ptheight0.5ptdepth0pt}\vskip2pt}\hskip.5pt}
\def\vt{{\tau\hskip-4.55pt\iota\hskip.6pt}}
\def\hs{\hskip.7pt}
\def\hh{\hskip.4pt}
\def\nh{\hskip-.7pt}
\def\hn{\hskip-.4pt}
\def\w{^{\phantom i}}

\def\bna{\hs\overline{\nh\nabla\nh}\hs}

\def\bg{\hskip1.2pt\overline{\hskip-1.2ptg\hskip-.3pt}\hskip.3pt}

\def\brc{\hskip.2pt\overline{\hskip-.2pt\mathrm{r}\nh}\hs}

\def\bsc{\hskip.2pt\overline{\hskip-.2pt\mathrm{s}\nh}\hs}

\def\eu{e_1\w}
\def\ed{e_2\w}
\def\et{e_3\w}
\def\eq{e\hn_4\w}

\newtheorem{thm}{\bf Theorem}[section]
\newtheorem{lem}[thm]{\bf Lemma}

\theoremstyle{definition}

\theoremstyle{remark}
\newtheorem{rem}[thm]{\bf Remark}

\numberwithin{equation}{section}

\begin{document}
\begin{abstract}
We call a Riemannian metric on a four-man\-i\-fold, as well as the underlying 
con\-for\-mal structure, {\it con\-for\-mal\-ly even}\hh$/\nnh${\it odd\/} if
its self-dual and anti-self-dual Weyl tensors acting on bi\-vec\-tors have,
at each point, equal/op\-po\-site spectra. The term {\it mixed parity\/}
refers to the simultaneous occurrence of both of these conditions, an
obvious example being the con\-for\-mal\-ly flat case. 
After listing various known classes of con\-for\-mal\-ly even metrics, we
construct examples of 
con\-for\-mal\-ly odd ones which are not of mixed parity, as well as some 
metrics of mixed parity that are not con\-for\-mal\-ly flat. The former 
examples are provided by left-in\-var\-i\-ant metrics on Lie groups,  
and realize uncountably many local con\-for\-mal-equiv\-a\-lence types.  
The latter also include some arising from Lie groups, as well as 
certain $1+3$ Riemannian products that are not
locally con\-for\-mal\-ly homogeneous.
\end{abstract}

\maketitle

\section{Introduction}
A con\-for\-mal structure on a four-di\-men\-sion\-al manifold, as well as
any Riemannian metric within it, may be called
{\it con\-for\-mal\-ly even\/} (or, {\it odd\/}) if its self-dual and
anti-self-dual Weyl tensors acting on bi\-vec\-tors have equal (or,
respectively, op\-po\-site) spectra at each point. We further speak of
{\it mixed parity\/} when both of these conditions occur at once, as they do,
for instance, in the con\-for\-mal\-ly flat case.

An obvious motivation behind these concepts is the natural question about how
changing the orientation may affect a given con\-for\-mal structure.

Our terminology is unrelated to the
occasional uses of the term `con\-for\-mal\-ly even'
when referring to polynomials in algebra and number theory, or in the context 
of asymptotically even con\-for\-mal\-ly compact metrics.

Here are some well-known classes -- see Section~\ref{ac} for details --
of con\-for\-mal\-ly even Riemannian
manifolds, the dimension being always assumed equal to four:
\begin{equation}\label{wkc}
\begin{array}{rl}
\mathrm{a)}&\mathrm{all\ con\-for\-mal\-ly\ flat\ manifolds,}\\
\mathrm{b)}&\mathrm{all\ warped\ products,\ including\ Riemannian\
products,}\\
\mathrm{c)}&\mathrm{all\ hypersurfaces\ in\ space\ forms,}\\
\mathrm{d)}&\mathrm{all\ manifolds\ admitting\ a\ Co\-daz\-zi\ tensor\ field,\
nowhere\ equal\ to}\\
&\mathrm{a\ multiple\ of\ the\ metric,}\\
\mathrm{e)}&\mathrm{all\ non}\hyp\mathrm{Ein\-stein\ manifolds\ having\
harmonic\ curvature,}\\
\mathrm{f)}&\mathrm{manifolds\ with\ pure\ curvature\ operator\ in\ the\
sense\ of\ Maillot\ }\text{\rm\cite{maillot}.} 
\end{array}
\end{equation}
In Section~\ref{co} we construct examples of Lie groups with 
left-in\-var\-i\-ant metrics which are con\-for\-mal\-ly odd
without being at the same time con\-for\-mal\-ly even. The Lie
al\-ge\-bras involved have Heis\-en\-berg-type derived ideals,
and the metrics form a one-pa\-ra\-me\-ter family of distinct 
local con\-for\-mal-equiv\-a\-lence classes.

Both Sections~\ref{co} and~\ref{mp} describe various mix\-ed-par\-i\-ty
metrics which are not con\-for\-mal\-ly flat, so that, at each point, 
$\,W^+\nnh$ and $\,W^-\nnh$ have the same spectrum
$\,\{-\sy,0,\sy\}\,$ with $\,\sy\ne0$. 
These examples also include some left-in\-var\-i\-ant metrics on Lie groups, 
in addition to certain $\,1+3\,$ Riemannian products 
without the lo\-cal-ho\-mo\-gene\-ity property,
even on the con\-for\-mal level.

\section{Preliminaries}\label{pr}
\setcounter{equation}{0}
Our sign convention about the curvature tensor $\,R\,$ of a Riemannian
metric $\,g$ is such that, for tangent vector fields $\,u,v,w\,$ and
the Le\-vi-Ci\-vi\-ta connection $\n\hs$, 
\begin{equation}\label{rvw}
R\hs(v,w)u\,=\,\nabla\!_{w}\w\n_v\w u\,-\,
\n_v\w\n_w\w u\,+\,\n_{[v,w]}\w u\hh.
\end{equation}
Given an algebraic curvature tensor $\,R\,$ in 
a Euclidean space $\,\mathcal{T}\hs$ of any dimension $\,n\ge4$, denoting by
$\,g\,$ the inner product, by $\,W\nnh\nh,\ric\,$ the Weyl and 
Ric\-ci tensors of $\,R$, and by $\,\sca\,$ its scalar curvature, we have
\begin{equation}\label{wey}
\begin{array}{l}
W_{\!ijpq}\w=R_{ijpq}\w
-\displaystyle{\frac1{n-2}}(g_{ip}\w\ric\nnh_{jq}\w
+g_{jq}\w\ric\nnh_{ip}\w-g_{jp}\w\ric\nnh_{iq}\w
-g_{iq}\w\ric\nnh_{jp}\w)\\
\hskip61.5pt
+\,\displaystyle{\frac{\sca}{(n-1)(n-2)}}
(g_{ip}\w g_{jq}\w-
g_{jp}\w g_{iq}\w).\phantom{1^{1^1}}
\end{array}
\end{equation}
If $\,n=4$,
since $\,W\hs$ is well known \cite[p.\,647]{hdg00} to commute with the 
Hodge star, the components of $\,W\nh$ in any or\-tho\-nor\-mal basis 
$\,\eu,\ed,\et,\eq$ satisfy the condition
\begin{equation}\label{wij}
\begin{array}{l}
W_{\!ijkl}\w\,=\,W_{\!pqrs}\w\mathrm{\ whenever\ both\
}\,\hs(i,j,p,q)\hs\,\mathrm{\ and}\\
(k,l,r,s)\mathrm{\ are\ even\ permutations\ of\ }(1,2,3,4)\hn.
\end{array}
\end{equation}
For the Ric\-ci tensors $\,\brc,\hs\ric$, 
scalar curvatures $\,\bsc,\hs\sca$, and La\-place operators
$\,\bde\hs,\hh\Delta\,$ (with pos\-i\-tive-def\-i\-nite
symbols) of two con\-for\-mal\-ly related 
metrics $\,\bg\,$ and $\,g=\bg/\hi^2$ in
dimension $\,n$, 
and any function $\,\psi$, one has the well-known relations
\begin{equation}\label{rct}
\begin{array}{rl}
\mathrm{i)}&\ric=\brc+(n-2)\hi^{-\nh1}\hs\bna\nh d\hi
+[\hi^{-\nh1}\bde\hi
-(n-1)\hi^{-\nh2}\hs\bg(\bna\nnh\hi,\nnh\bna\nnh\hi)]{\bg}\hh,\\
\mathrm{ii)}&\sca\,=\,\hi^2\hs\bsc+2(n-1)\hi\hs\bde\hs\hi
-n(n-1)\hs\bg(\bna\nnh\hi,\nnh\bna\nnh\hi)\hh.\\
\mathrm{iii)}&\Delta\psi=\hi^2\hh\bde\hh\psi
-(n-2)\hi\hs\bg(\bna\nh\hi,\nnh\bna\nh\psi)\hh,\,\,\,
g(\nabla\nh\psi,\nnh\nabla\nh\psi)
=\hi^2\nh\bg(\bna\nnh\psi,\nh\bna\nnh\psi).
\end{array}
\end{equation}
See, e.g., \cite[p.\,529]{hdg00} and
\cite[formulae (7.1), (7.2), (7.6)]{derdzinski-park-shin-26}. Here
$\,\bna\nnh$ denotes both
the Le\-vi-Ci\-vi\-ta connection of $\,\bg\,$ and the $\,\bg\hh$-grad\-i\-ent. 
When $\,\hat g=\psi g$, (\ref{rct}-ii) becomes
\begin{equation}\label{scu}
4\psi^3\hat\sca\,=\,4\psi^2\hs\sca\,-\,(n-1)[4\psi\hh\Delta\psi\,
+\,(n-6)g(\nabla\nh\psi,\nabla\nh\psi)].
\end{equation}
\begin{lem}\label{lineq}
With\/ $\,i,j,p,q\in\{1,2,3,4\}\,$ always assumed distinct, and any six
given real scalars\/ $\,\sy\nnh_{ij}\w$, $\,i<j$, the system
$\,2(\ky_i\w+\ky_j\w)-(\ky_p\w+\ky_q\w)=6\hh\sy\nnh_{ij}\w$ of six linear
equations with four unknowns\/ $\,\ky_1\w,\dots,\ky_4\w$ has a solution if and only if
\begin{equation}\label{iff}
\sy\nnh_{12}\w\hn+\sy\nnh_{34}\w\,=\,\sy\nnh_{13}\w\hn+\sy\nnh_{24}\w\,
=\,\sy\nnh_{14}\w\nh+\sy\nnh_{23}\w,
\end{equation}
and then the solution is unique, namely, 
$\,\ky_i\w=\sy\nnh_{ij}\w+\sy\nnh_{ip}\w+\sy\nnh_{iq}\w$.  
\end{lem}
\begin{proof}[{\it Proof}]Uniqueness is obvious if we use the system to
evaluate $\,\sy\nnh_{ij}\w+\sy\nnh_{ik}\w+\sy\nnh_{ip}\w$,
which also shows that the
$\,6\times4\,$ matrix of the system has rank four. The system further implies
that $\,6(\sy\nnh_{ij}\w\hn+\sy\nnh_{pq}\w)\,$ are all equal to
$\,\ky_1\w+\ldots+\ky_4\w$,
proving the `only if' claim. The `if' part follows: the image of our
rank-four matrix operator $\,\bbR\hn^4\nh\to\bbR\hn^6\nh$, being contained
in the sub\-space defined by the two independent equations (\ref{iff}), must
coincide with it. This completes the proof.
\end{proof}
\begin{rem}\label{cfhom}We say that a Riemannian manifold
$\,(M\nh,g)\,$ of dimension $\,n\,$ is {\it locally con\-for\-mal\-ly
homogeneous\/} when any two points $\,x,x'\in M\,$ have neighborhoods
$\,\,U\nh,\,U'$ with a con\-for\-mal dif\-feo\-mor\-phism 
$\,\,U\nh\to U'$ sending $\,x\,$ to $\,x'\nh$. If this is the case and 
$\,n\ge4\,$ (or, $\,n=4\,$ and $\,M\,$ is oriented), while 
$\,W\nnh\nh\ne\hs0\,$ (or, $\,W^\pm\nnh\nh\ne\hs0$) everywhere,
then the metric $\,|W\hn|\hh g\,$ (or, $\,|W^\pm\nnh|\hh g$) is
locally homogeneous.
In fact, the new metric is con\-for\-mal\-ly invariant.
\end{rem}

\section{Algebraic curvature tensors and con\-for\-mal evenness}\label{ac}
\setcounter{equation}{0}
For an algebraic curvature tensor $\,R\,$ in an oriented Euclidean
four-space $\,\mathcal{T}\nnh\nh$, we always denote by 
$\,W\hskip-3pt,\,\sca\,$ and $\,\ric\,$ the Weyl component, scalar curvature,
and Ric\-ci contraction of $\,R$, by $\,\varLambda\nnh^\pm$ the spaces of
self-dual and anti-self-dual bi\-vec\-tors in $\,\mathcal{T}\nnh\nh$, and
by $\,W^\pm\nnh$ the self-dual and anti-self-dual parts of $\,W\nnh\nh$.
Any positive orthonormal basis $\,\eu,\ed,\et,\eq$ of 
$\,\mathcal{T}\hs$ gives rise to the pair
\begin{equation}\label{bss}
\pm\hh\eu\wedge\ed+\et\wedge\eq\hh,\quad 
\pm\hh\eu\wedge\et+\eq\wedge\ed\hh,\quad 
\pm\hh\eu\wedge\eq+\ed\wedge\et,
\end{equation}
of length $\,\sqrt{2\,}$ orthogonal bases
of $\,\varLambda\nnh^\pm\nnh$, thus making both $\,\varLambda\nnh^\pm$ 
canonically oriented. Conversely, it is well known -- see, for instance,
\cite[Lemma 7.1]{derdzinski-park-shin-27} -- that
\begin{equation}\label{evp}
\begin{array}{l}
\mathrm{every\ pair\ of\ length\ }\hs\sqrt{2\hh}\mathrm{\ positive\ orthogonal\
bases\ of\ the\ spaces\ }\,\varLambda\nnh^\pm\\ 
\mathrm{equals\hh\ (\ref{bss})\hh\ for\hh\ some\hh\ positive\hh\ 
or\-tho\-nor\-mal\hh\ basis\hh\ }\,\eu,\ed,\et,\eq\mathrm{\hh\ of\hh\
}\,\mathcal{T}\nnh.
\end{array}
\end{equation}
\begin{lem}\label{diagw}
For\/ $\,R,W^\pm\nnh$ and\/ $\,\eu,\ed,\et,\eq$ as above, 
the bases\/ {\rm(\ref{bss})} consist of eigen\-vectors of\/ $\hs W^\pm\nnh$ 
acting on bi\-vec\-tors if and only if, with\/
$\,W_{\!ijpq}\w=W\hn(e\nh_i\w,e\nh_j\w,e\hn_p\w,e\hn_q\w)$, 
\begin{equation}\label{thr}
W_{\!ijpq}\w=0\,\mathrm{\ whenever\ }\,\{i,j,p,q\}\,\mathrm{\ has\ three\
elements,}
\end{equation}
and then the respective eigen\-values are\/ 
$W_{\!1212}\w\nh\pm\hn W_{\!1234}\w,\hs W_{\!1313}\w\nh\pm\hn W_{\!1342}\w,\hs
W_{\!1414}\w\nh\pm\hn W_{\!1423}\w$.
\end{lem}
\begin{proof}[{\it Proof}]
If (\ref{bss}) are eigen\-vectors of $\,W^\pm\nnh$, each
pair formed by $\,\eu\wedge e\hn_i\w+e\nh_j\w\wedge e\hn_k\w$ and 
$\,-\eu\wedge e\hn_i\w+e\nh_j\w\wedge e\hn_k\w$, where 
$\,(i,j,k)\,$ is any cyclic permutation of $\,(2,3,4)$, spans
a $\,W\hn\nnh$-in\-var\-i\-ant space of bi\-vectors, and hence so does the
pair $\,\eu\wedge e\hn_i\w,\,\,e\nh_j\w\wedge e\hn_k\w$, leading to 
(\ref{thr}). As (\ref{thr}) implies the $\,W\hn\nnh$-in\-var\-i\-ance of each
such span of $\,\eu\wedge e\hn_i\w,\,\,e\nh_j\w\wedge e\hn_k\w$, that is, 
of the span of $\,\eu\wedge e\hn_i\w+e\nh_j\w\wedge e\hn_k\w$ and 
$\,-\eu\wedge e\hn_i\w+e\nh_j\w\wedge e\hn_k\w$, the converse follows 
from the fact that the sub\-spaces $\,\varLambda\nnh^\pm$ are also 
$\,W\nnh$-in\-var\-i\-ant, combined with (\ref{wij}).
\end{proof}
Following \cite[Sect.~16.18]{besse}, one says that an
or\-tho\-nor\-mal basis $\,\eu,\ed,\et,\eq$ of $\,\mathcal{T}$ {\it
di\-ag\-o\-nal\-izes\/} an algebraic curvature tensor $\,R\,$ if 
$\,R(e\nh_i\w,e\nh_j\w,e\hn_p\w,e\hn_q\w)=0$ whenever 
$\,\{i,j\}\ne\{p,q\}\,$ or, in other words, if all the
nonzero exterior products $\,e\hn_i\w\nh\wedge e\nh_j\w$ are eigenvectors of 
$\,R\,$ acting on bi\-vec\-tors. Riemannian 
manifolds $\,(M\nh,g)\,$ of all dimensions, whose curvature tensor is
analogously di\-ag\-o\-nal\-iz\-ed, at every point $\,x$, by some
or\-tho\-nor\-mal basis of $\,\txm\nh$, were first studied by Maillot
\cite{maillot}, who referred to them as having {\it pure curvature 
operator}. 

We call $\,R\,$ {\it con\-for\-mal\-ly even\/} when the self-dual and
anti-self-dual parts $\,W^\pm\nnh$ of its Weyl component, 
acting on $\,\varLambda\nnh^+$ and $\,\varLambda\nnh^-\nnh\nh$,
have the same spectrum.
\begin{lem}\label{cfevn}
An algebraic curvature tensor\/ $\,R\,$ in a Euclidean four-space\/
$\,\mathcal{T}\hs$ is con\-for\-mal\-ly even if and only if some
or\-tho\-nor\-mal basis\/ $\,\eu,\ed,\et,\eq$ of\/ $\,\mathcal{T}\hs$ 
di\-ag\-o\-nal\-izes its Weyl component\/ $\,W\nnh\nh$.
\end{lem}
\begin{proof}[{\it Proof}]
The `if' part is \cite[Lemma 16.20]{besse}. It also 
follows from 
Lemma~\ref{diagw}. On the other hand, when $\,R\,$ is con\-for\-mal\-ly even,
using (\ref{evp}) to select an or\-tho\-nor\-mal basis $\,\eu,\ed,\et,\eq$ of
$\,\mathcal{T}\hs$
for which (\ref{bss}) are eigen\-vectors of $\,W^\pm\nnh\nh$, with
the eigenvalues listed in Lemma~\ref{diagw} that now do not depend on the sign
$\,\pm\hs$, we get (\ref{thr}) and 
$\,W_{\!1234}\w=W_{\!1342}\w=W_{\!1423}\w$, as required.
\end{proof}
Our claim about (\ref{wkc}-a) is obvious. As for (\ref{wkc}-b) in the case 
of $\,1+3\,$ Riemannian products, con\-for\-mal evenness
clearly follows when every point $\,x\,$ is a fixed point of some
o\-ri\-en\-ta\-tion-re\-vers\-ing isometry of a neighborhood of $\,x$.

In a K\"ah\-ler surface with the standard o\-ri\-en\-ta\-tion and the scalar
curvature $\,\sca$, the spectrum of $\,W^+\nnh$ is
$\,\{\sca/6,-\sca/\nh12,-\sca/\nh12\}$, 
the eigen\-val\-ue function $\,\sca/6\,$ being realized by the K\"ah\-ler form
treated as a bi\-vec\-tor. See, for instance., \cite[Proposition 9.8]{hdg00}.
This yields (\ref{wkc}-b) 
for a Riemannian product of two surfaces, as it is locally 
K\"ah\-ler for two complex structures inducing opposite orientations.
The case of warped products is now immediate from con\-for\-mal invariance.

An or\-tho\-nor\-mal basis di\-ag\-o\-nal\-iz\-ing $\,R\,$ also
di\-ag\-o\-nal\-izes its Weyl component and -- in the usual sense --
its Ric\-ci tensor. Lemma~\ref{cfevn} now proves our claim about
(\ref{wkc}-f) and -- due to the Gauss equation -- also (\ref{wkc}-c),
while for (\ref{wkc}-d) it suffices to invoke \cite[p.\,440,\,16.21]{besse}
and, again, Lemma~\ref{cfevn}. Finally, combining (\ref{wkc}-d) and
(\ref{wkc}-b) with \cite[p.\,432, Chapter\,16\hs A]{besse}, we obtain the
claim about (\ref{wkc}-e).

\section{Con\-for\-mal\-ly odd curvature tensors}\label{oc}
\setcounter{equation}{0}
Given an algebraic curvature tensor $\,R\,$ in an oriented Euclidean 
four-space $\,\mathcal{T}\nnh\nh$, we say that $\,R\,$ is
{\it con\-for\-mal\-ly odd\/} if its self-dual and anti-self-dual Weyl
components $\,W^\pm\nnh$ acting on $\,\varLambda\nnh^+$ and 
$\,\varLambda\nnh^-$ (see Section~\ref{ac}) have op\-po\-site spectra.
\begin{lem}\label{cfodd}
Any\/ $\,R\,$ as above is con\-for\-mal\-ly odd if and only if,
for some positive orthonormal basis\/ $\,\eu,\ed,\et,\eq$ of\/  
$\,\mathcal{T}\nnh\nh$, whenever\/ $\,\{i,j,p,q\}=\{1,2,3,4\}$,
\begin{equation}\label{cnd}
\mathrm{a)}\hskip5ptR\hn_{ijij}\w\nh+\nh R\hn_{pqpq}\w\mathrm{\hs\ are\ all\
equal\ to\ one\ another,\ and\ \ }\mathrm{b)}\hskip5ptR\hn_{ipiq}\w\nh
=\nh R\hn_{jpjq}\w.
\end{equation}
As a condition imposed on the basis\/ $\,\eu,\ed,\et,\eq$, {\rm(\ref{cnd})}
amounts to requiring that the bases\/
{\rm(\ref{bss})} di\-ag\-o\-nal\-ize\/ $\,W^\pm\nnh\nnh$, while
\begin{equation}\label{ope}
\eu\wedge e\hn_i\w+e\nh_j\w\wedge e\hn_k\w\mathrm{\ and\ 
}\nnh-\nh\nnh\eu\wedge e\hn_i\w+e\nh_j\w\wedge e\hn_k\w\mathrm{\
correspond\ to\
opposite\ eigen\-values}
\end{equation}
of\/ $\,W^+\nnh$ and\/ $\,W^-\nnh\nnh$, for any cyclic permutation\/ 
$\,(i,j,k)\,$ of\/ $\,(2,3,4)$.

The spectra of\/ $\,W^\pm\nh$ realized by {\rm(\ref{ope})} then are\/ 
$\,\pm R\hn_{1234}\w,\,\pm R\hn_{1342}\w,\,\pm R\hn_{1423}\w$.
\end{lem}
\begin{proof}[{\it Proof}]
Let $\,R\,$ be con\-for\-mal\-ly odd. Using (\ref{evp}) we may 
choose $\,\eu,\ed,\et,\eq$ for which (\ref{bss}) di\-ag\-o\-nal\-izes
$\hh W^\pm\nnh\hn$ so that (\ref{ope}) holds.
Then Lemma~\ref{diagw} and (\ref{wij}) give $\,W_{\!ijpq}\w=0\,$ if 
$\,\{i,j,p,q\}\ne\{1,2,3,4\}$. In view of (\ref{wey}), this
amounts to
\begin{equation}\label{amt}
\mathrm{a)}\hskip5pt6R\hn_{ijij}\w\nh
=\nh2(\ric\nnh_{ii}\w\nh+\nh\ric\nnh_{jj}\w)\nh
-(\ric\nnh_{pp}\w\nh+\nh\ric\nnh_{qq}\w)\mathrm{\ \ \ and\ \ \ 
}\mathrm{b)}\hskip5pt2R\hn_{ipiq}\w\nh=\nh\ric\nnh_{pq}\w
\end{equation}
whenever $\,\{i,j,p,q\}=\{1,2,3,4\}$, which is equivalent to (\ref{cnd}),
since (\ref{cnd}-a) implies (\ref{amt}-a) in view of
Lemma~\ref{lineq} applied to $\,\sy\nnh_{ij}\w=R\hn_{ijij}\w$ and
$\,\ky_i\w=\ric\nnh_{ii}\w$.
 
Conversely, assuming (\ref{cnd}), we get, as stated above, (\ref{amt}) and,
therefore, $\nh W_{\!ijpq}\w\nnh=0\,$ if 
$\,\{i,j,p,q\}\ne\{1,2,3,4\}$. By Lemma~\ref{diagw}, the basis in question
consist of eigen\-vectors of $\hs W^\pm\nnh$ with the eigen\-values
$\,\pm W_{\!1234}\w,\,\pm W_{\!1342}\w,\,\pm W_{\!1423}\w$, as required.
\end{proof}

\section{Con\-for\-mal\-ly odd metrics}\label{co}
\setcounter{equation}{0}
We now construct examples of metrics that are
con\-for\-mal\-ly odd without being con\-for\-mal\-ly even, as well as
some of mixed parity that are not con\-for\-mal\-ly flat.

We treat a left-in\-var\-i\-ant metric $\,g\,$ 
on a Lie group as an inner product in the Lie algebra
of left-in\-var\-i\-ant vector fields, and always introduce  
the Lie-al\-ge\-bra structure by providing the brackets of
vectors from a fixed or\-tho\-nor\-mal basis.
\begin{thm}\label{xcfod}
Given\/ $\,\ky,\ly,\my\in\bbR\,$ with\/ $\,2\ky^2\nh=-\ly\my$,
the brackets defined,
in an or\-tho\-nor\-mal basis\/ $\,\x,\y,\z,\t$, by
\begin{equation}\label{brd}
[\x,\y]=\ly\y,\,\,[\x,\z]=\my\nh\z,\,\,
{}[\x,\t]=(\ly+\my)\t,\,\,[\y,\z]=2\ky\hh\t,\,\,[\y,\t]=[\z,\t]=0,
\end{equation}
satisfy the Ja\-co\-bi
identity. For the left-in\-var\-i\-ant metric\/ $\,g\,$ 
on a four-di\-men\-sion\-al Lie group, arising as stated above,
and the orientation determined by\/ $\,\x,\y,\z,\t$, 
the spectra of\/ $\,W^\pm\nnh$ are\/ 
$\,\{\pm\ly\ky,\,\pm\my\ky,\,\mp(\ly+\my)\ky\}$. Furthermore,
\begin{enumerate}
\item[{\rm(i)}] $g\,$ is con\-for\-mal\-ly odd, but not con\-for\-mal\-ly
even, when\/ $\,(\ly+\my)\ly\my\ne0$,
\item[{\rm(ii)}] $g\,$ has mixed parity without being con\-for\-mal\-ly flat
if\/ $\,\my=-\nnh\ly\ne0$,
\item[{\rm(iii)}] $g\,$ is con\-for\-mal\-ly flat when\/ $\,\ky=0\,$ or\/
$\,\ly=\my=0$.
\end{enumerate}
\end{thm}
\begin{proof}[{\it Proof}]
The Ja\-co\-bi identity is obvious, and 
the Le\-vi-Ci\-vi\-ta connection $\,\n\,$ has
\begin{equation}\label{lcc}
\begin{array}{l}
\n_\x\w\x=\n_\x\w\y=\n_\x\w\z=\n_\x\w\t=0,\,\,\,\n_\y\w\x=-\ly\y,\,\,\,\n_\y\y
=\ly\x,\\
\n_\y\w\z=\ky\hh\t,\quad\n_\y\w\t=-\ky\z,\quad\n_\z\w\x
=-\my\z,\quad\n_\z\w\y=-\ky\hh\t,\\
\n_\z\w\z=\my\x,\qquad\n_\z\w\t=\ky\y,\qquad\n\hs_\t\w\x=-(\ly+\my)\t,\\
\n\hs_\t\w\y=-\ky\z,\qquad\n\hs_\t\w\z=\ky\y,\qquad\n\hs_\t\w\t=(\ly+\my)\x
\end{array}
\end{equation}
since, for $\,\n\,$ {\it defined\/} by (\ref{lcc}) and all $\,u,v\,$ from
our basis, $\,\n_u\w$ is skew-adjoint and $\,\n_u\w v-\n_v\w u=[u,v]$.
Using (\ref{rvw}), we easily verify that, with
$\,R_{\x\y\x}\w=R(\x,\y,\x)$, etc.,
\[
\begin{array}{l}
R_{\x\y\x}\w=-\ly^2\nh\y,\,\,\,
R_{\x\y\y}\w=\ly^2\nh\x,\,\,\,
R_{\x\y\z}\w=\ly\ky\hh\t,\,\,\,
R_{\x\y\t}\w=-\ly\ky\z,\\
R_{\x\z\x}\w=-\my^2\nh\z,\,\,\,
R_{\x\z\y}\w=-\my\ky\hh\t,\,\,\,
R_{\x\z\z}\w=\my^2\nh\x,\,\,\,
R_{\x\z\t}\w=\my\ky\y,\\
R_{\x\t\x}\w=-(\ly+\my)^2\hn\t,\,\,
R_{\x\t\y}\w=-(\ly+\my)\ky\z,\,\,
R_{\x\t\z}\w=(\ly+\my)\ky\y,\\
R_{\x\t\t}\w=(\ly+\my)^2\nh\x,\,\,
R_{\y\z\x}\w=-(\ly+\my)\ky\hh\t,\,\,
R_{\y\z\y}\w=-(3\ky^2\nnh+\ly\my)\z,\\
R_{\y\z\z}\w=(3\ky^2\nnh+\ly\my)\y,\qquad
R_{\y\z\t}\w=(\ly+\my)\ky\x,\\
R_{\y\t\x}\w=-\my\ky\z,\quad
R_{\y\t\y}\w=(\ky^2\nnh-\ly^2\nnh-\ly\my)\t,\quad
R_{\y\t\z}\w=\my\ky\x,\\
R_{\y\t\t}\w=(\ly^2\nnh+\ly\my-\ky^2)\y,\quad
R_{\z\t\x}\w=\ly\ky\y,\quad
R_{\z\t\y}\w=-\ly\ky\x,\\
R_{\z\t\z}\w=(\ky^2\nnh-\ly\my-\my^2)\t,\qquad
R_{\z\t\t}\w=(\ly\my+\my^2\nnh-\ky^2)\z.
\end{array}
\]
Therefore, (\ref{cnd}-b) is satisfied due to
the six obvious equalities
\[
\begin{array}{l}
R_{\x\z\x\t}\w=R_{\y\z\y\t}\w=0,\qquad
R_{\x\y\x\t}\w=R_{\z\y\z\t}\w=0,\qquad
R_{\x\y\x\z}\w=R_{\t\y\t\z}\w=0,\\
R_{\y\x\y\t}\w=R_{\z\x\z\t}\w=0,\qquad
R_{\y\x\y\z}\w=R_{\t\x\t\z}\w=0,\qquad
R_{\z\x\z\y}\w=R_{\t\x\t\y}\w=0.
\end{array}
\]
Finally, (\ref{cnd}-a) reads 
$\,\ky^2\nnh-\ly\my-\my^2\hs=\,\ky^2\nnh-\ly\my-\my^2\hs=\,-3\ky^2\nnh
-3\ly\my-\my^2\nh$, 
which is nothing else than $\,2\ky^2\nh=-\ly\my$.
Lemma~\ref{cfodd} now yields our assertion
about the spectra of $\,W^\pm\nnh\nh$, and (i) --- (iii) trivially
follow.
\end{proof}
The examples arising from part (i) of Theorem~\ref{xcfod}
represent {\it an uncountable family of local con\-for\-mal types\/}: 
generally, the ratio of two nonzero eigen\-values of $\,W^\pm\nh$ is 
a con\-for\-mal invariant, and our $\,\my/\nh\ly\,$ assume all
values in $\,(-\infty,-\nnh1)\cup(-\nnh1,0)$.

\section{Further examples of mixed parity}\label{mp} 
\setcounter{equation}{0}
As pointed out in (\ref{wkc}-b), all Riemannian products 
of dimension four are con\-for\-mal\-ly even. Some such products, of the form
\begin{equation}\label{prd}
(\bbR\times\nh M\nh,\,\,dt^2\nnh+\,g),
\end{equation}
are con\-for\-mal\-ly odd as well (mixed parity) -- both
$\,W^\pm\nnh$ have the same unordered spectrum 
$\,\{-\sy,0,\sy\}$, with $\,\sy\ne0\hh$:
\begin{lem}\label{mixed}
In the Riemannian product of a three-man\-i\-fold\/ $\,(M\nh,g)\,$ and a
one-di\-men\-sion\-al factor, the self-dual and anti-self-dual Weyl tensors\/
$\,\,W^\pm\nnh\nnh$ acting on bi\-vec\-tors are nonzero and
have op\-po\-site spectra at precisely those points where
\begin{equation}\label{ric}
\begin{array}{l}
\mathrm{one\hh\ eigen\-value\hh\ of\hh\ the\hh\ Ric\-ci\hh\ tensor\hh\ of\hh\ 
}\,(M\nh,g)\,\mathrm{\hh\ 
is\hh\ the\hh\ arith}\hyp\\
\mathrm{metic\nh\ mean\nh\ of\nh\ the\nh\ other\nh\ two,\nnh\nh\ without\nh\
being\nh\ equal\ to\nh\ them.}
\end{array}
\end{equation}
\end{lem}
In fact, one easily verifies -- see, e.g., 
\cite[formula (8.4) with $\varphi=1$]{derdzinski-park-shin-26} -- that  
\begin{equation}\label{spc}
\begin{array}{l}
\mathrm{at\ each\ point\ both}-\hskip-1.7pt2W^\pm\nnh\mathrm{\ have\ the\
same\ eigen\-values,}\\
\mathrm{equal\ to\ those\ of\ the\ trace\-less\ Ric\-ci\ tensor\ of\
}\,(M\nh,g).
\end{array}
\end{equation}
Due to Lemma~\ref{mixed}, {\it three-man\-i\-folds\/
$\,(M\nh,g)\,$ with\/ {\rm(\ref{ric})} at every point 
lead to examples of the form\/ {\rm(\ref{prd})}
that have mixed parity without being con\-for\-mal\-ly flat.}

One class of such $\,(M\nh,g)\,$ is provided by certain Lie groups of
dimension three with left-in\-var\-i\-ant metrics, introduced as
in the lines preceding Theorem~\ref{xcfod}: the brackets are given,
in an or\-tho\-nor\-mal basis\/ $\,\y,\z,\t\,$ of the Lie algebra, by
\begin{equation}\label{yze}
[\y,\z]=2\hh\yc\hs\t,\qquad[\z,\t]=2\ya\y,\qquad[\t,\y]=2\hh\yb\z,
\end{equation}
for any fixed $\,\ya,\yb,\yc\in\bbR$, and then, according to
\cite[Theorem 4.3]{milnor}, 
the basis $\,\y,\z,\t\,$ di\-ag\-o\-nal\-izes the Ric\-ci tensor, with
the eigen\-values $\,2\beta\gamma,\,2\alpha\gamma,\,2\alpha\beta$,
where $\,(\alpha,\beta,\gamma)=(\ya-\yb-\yc,\,\yb-\ya-\yc,\,\yc-\ya-\yb)$.
To obtain (\ref{ric}), we just note that 
the image of the Cre\-mo\-na transformation
sending $\,(\alpha,\beta,\gamma)\in\rtr$ to
$\,(\beta\gamma,\alpha\gamma,\alpha\beta)$ obviously consists of triples 
having three positive entries, or one positive and two negative entries,
or at most one nonzero entry.

According to the final clause of Theorem~\ref{nchom} below, 
the other class fundamentally differs both from the preceding one and 
from that arising in Theorem~\ref{xcfod}(ii).

Namely, in our second class of examples, $\,(M\nh,g)\,$ is con\-for\-mal\-ly
flat. We fix a Euclidean three-space $\,V\nnh$ with the inner product
$\,\langle\,,\rangle\,$ and a self-ad\-joint linear en\-do\-mor\-phism
$\,\hj\nh$, which is not a multiple of $\,\mathrm{Id}$, such that 
\begin{equation}\label{one}
\mathrm{one\ 
eigen\-value\ of\ }\,\hj\,\mathrm{\ equals\ the\
arithmetic\ mean\ of\ the\ other\ two.}
\end{equation}
Denoting by $\,\bg\,$ the flat metric represented by
$\,\langle\,,\rangle$, we define $\,g\,$ by $\,g=\bg/\hi^2\nh$,
where $\,\hi:V\to\bbR$, with $\,\hi(v)=\langle\hj\hn v,v\rangle$, on 
any connected component $\,M\,$ of 
$\,V\nh\smallsetminus\hi^{-\nnh1}(0)$.
\begin{thm}\label{nchom}
If\/ $\,\mathrm{tr}\hskip1pt\hj\ne0$, that is,
$\,\hj\hh$ has a nonzero middle eigen\-value, then\/ {\rm(\ref{prd})} for\/
$\,(M\nh,g)\,$ defined above has mixed parity without being con\-for\-mal\-ly
flat, and is not locally 
con\-for\-mal\-ly homogeneous, in the sense of Remark\/~{\rm\ref{cfhom}}.
\end{thm}
\begin{proof}[{\it Proof}]
For $\,\vt=\mathrm{tr}\hskip1pt\hj\,$ and $\,Q:M\to\bbR\,$ given
by $\,Q(v)=\langle\hj\hn v,\hj\hn v\rangle$, we have (\ref{rct}-i) with 
$\,(n,\,\brc,\,\bna\nh d\hi,\,\bde\hs\hi,\bg(\bna\nnh\hi,\nnh\bna\nnh\hi))\,
=\,(3,\,0,\,2\langle\hj\hs\cdot\,,\,\cdot\,\rangle,\,2\vt,\,4Q)$, 
so that (\ref{one}) yields
(\ref{ric}), as the conclusion is obviously
unaffected by the term proportional to $\,\bg\,$ in (\ref{rct}-i). Also,
from (\ref{rct}-ii) -- (\ref{rct}-iii), on $\,(M\nh,g)$,
\begin{equation}\label{scc}
\sca\,=\,8(\vt\hh\hi\,-\,3Q),\quad\Delta\hi\,=\,2(\vt\hh\hi-2Q)\hh\hi,\quad
g(\nabla\nh\hi,\nnh\nabla\nh\hi)\,=\,4Q\hi^2\nh.
\end{equation}
By (\ref{rct}-i), the bundle morphism $\,T\nh M\to T\nh M\,$
corresponding, via $\,g$, to the trace\-less Ric\-ci tensor of $\,(M\nh,g)$,
equals $\,2\hi\,$ times the trace\-less part of $\,\hj\nh$. (We need not
invoke (\ref{rct}-ii), since the term proportional to
$\,\bg\,$ in (\ref{rct}-i) has zero trace\-less part.) Therefore, 
according to (\ref{spc}), $\,W^\pm\nnh$ both have the same spectrum 
$\,\{-\sy,0,\sy\}\,$ with $\,\sy=(\ly-\my)\hi/2$, where
$\,\ly,\my,(\ly+\my)/2\,$ are the eigen\-values of $\,\hj$. 
As $\,\hj\,$ is not a multiple of $\,\mathrm{Id}$, we have 
$\,\sy\ne0$, which proves the first part of the theorem. Also,
\begin{equation}\label{nrm}
\mathrm{the\ norm\ of\ }\,W^\pm\nnh\mathrm{\ equals\ a\ constant\ times\ 
}\,|\hi|.
\end{equation}
Combining (\ref{scc}), which obviously holds on the product manifold
(\ref{prd}) as well, with (\ref{scu}) for $\,(n,\psi)=(4,|\hi|)\,$ and $\,g\,$
replaced by $\,dt^2\nnh\nh+g$, we easily verify that the metric
$\,\hat g=|\hi|\hh(dt^2\nnh\nh+g)\,$ on $\,\bbR\times\nh M\,$ has the scalar
curvature $\,\hat\sca\,$ equal to $\,2\hskip2.3pt\mathrm{sgn}\hskip1pt\hi$
times $\,\vt-3Q\hi^{-\nnh1}\nh$. Thus,
$\,\hat\sca\,$ is nonconstant: otherwise, $\,\hj^2$ would be a multiple of
$\,\hj\nh$, and the resulting equality $\,\hj^2\nh=c\hj\,$ would obviously
give $\,\ky^2\nh=c\ky\,$ for each eigen\-value $\,\ky\,$ of $\,\hj\nh$, so
that $\,\hj\,$ would have at most two distinct eigen\-values, contrary to
(\ref{one}) and the line preceding it.

Consequently, $\,(\bbR\times\nh M\nh,\,|\hi|\hh(dt^2\nnh\nh+g))\,$ is not
locally homogeneous. The final clause of the theorem now follows from
(\ref{nrm}) and Remark~\ref{cfhom}.
\end{proof}

\ \
%
%

%

\noindent{\bf Funding.}  
J.H.\,Park was supported by the National Research Foundation of Korea
(NRF) grant funded by the Korea government (MSIT) (RS-2024-00334956).
S.\,Kim was supported by Basic Science Research
Program through the National Research
Foundation of Korea (NRF) funded by the Ministry of Education (No.
RS-2019-NR040081).


\end{document}